\documentclass[12pt]{article}
\usepackage[a4paper, hmargin=2.2cm, vmargin=2.5cm]{geometry}
\usepackage{amsfonts,amsmath}
\usepackage[authoryear]{natbib} \bibpunct{(}{)}{,}{a}{,}{,}

\newcommand{\bA}{\mbox{\boldmath $A$}}
\newcommand{\bI}{\mbox{\boldmath $I$}}
\newcommand{\bX}{\mbox{\boldmath $X$}}
\newcommand{\bY}{\mbox{\boldmath $Y$}}
\newcommand{\bZ}{\mbox{\boldmath $Z$}}
\newcommand{\ba}{\mbox{\boldmath $a$}}
\newcommand{\bb}{\mbox{\boldmath $b$}}

\newcommand{\bt}{\mbox{\boldmath $t$}}

\newcommand{\by}{\mbox{\boldmath $y$}}
\newcommand{\bzero}{\mbox{\boldmath $0$}}
\newcommand{\bmu}{\mbox{\boldmath $\mu$}}

\newcommand{\bSigma}{\mbox{\boldmath $\Sigma$}}
\newcommand{\bOmega}{\mbox{\boldmath $\Omega$}}
\newcommand{\bdelta}{\mbox{\boldmath $\delta$}}

\newcommand{\blambda}{\mbox{\boldmath $\lambda$}}
\newcommand{\btheta}{\mbox{\boldmath $\theta$}}

\newcommand{\cov}{\mbox{cov}}
\newcommand{\var}{\mbox{var}}

\begin{document}

\title{On mean and/or variance mixtures of normal distributions}

\author{Sharon X. Lee$^{1}$, Geoffrey J. McLachlan$^{2, \star}$}
\date{}

\maketitle

\begin{flushleft}
$^1$School of Mathematical Sciences, University of Adelaide, 
Adelaide, South Australia, 5005, Australia.\\
$^2$Department of Mathematics, University of Queensland, 
Brisbane, Queensland, Australia, 4072, Australia.\\
$^\star$ E-mail: g.mclachlan@uq.edu.au
\end{flushleft}

\abstract{Parametric distributions are an important part 
of statistics. There is now a voluminous literature 
on different fascinating formulations of flexible distributions. 
We present a selective and brief overview of 
a small subset of these distributions, 
focusing on those that are obtained by 
scaling the mean and/or covariance matrix of 
the (multivariate) normal distribution 
with some scaling variable(s). 
Namely, we consider the families of mean mixture, 
variance mixture, and mean-variance mixture 
of normal distributions. Its basic properties, 
some notable special/limiting cases, 
and parameter estimation methods are also described.}

\section{Introduction}
\label{sec:intro}

The normal distribution plays a central role 
in statistical modelling and data analysis, 
but real data rarely follow this classical distribution. 
The quest for more flexible distributions has led to 
an ever growing development in the literature of 
parametric distribution. In the past two decades or so, 
intense interest has been in the area of skew 
or asymmetric distributions; see, for example, 
the book edited by \citet{B001}, 
the monograph by \citet{B025}, 
and the papers by \citet{J021,J017} and \citet{aa20} 
for recent accounts of the literature on skew distributions. 
Many of these formulations belong to 
the class of skew-symmetric distributions, 
which is a generalization of 
the classical skew normal (SN) distribution 
by \citet{J001}. 
This SN distribution can be characterized as 
a mean mixture of normal (MMN) distribution, 
where the mean of a normal random variable 
is scaled by a truncated normal random variable \citep{njsb19}. 
Another related and extensively studied family of distributions 
that can render asymmetric distributional shapes 
is the mean-variance mixture of normal (MVMN) distribution. 
Introduced by \citet{bns82}, the MVMM distribution 
is obtained by scaling both mean and variance 
of a normal random variable with 
the same (positive scalar) scaling random variable. 

This paper presents a brief overview of flexible distributions 
that arise from scaling either/both the mean 
and variance of a normal random variable. 
For simplicity, we focus on the case of 
a univariate scaling variable. 
Apart from the aforementioned MMN and MVMN families, 
a third family called variance mixture of normal 
(VMN) distributions can be defined by 
scaling only the variance of a normal random variable. 
Although VMN does not produce asymmetric distributions 
(at least not in the case of a scalar scaling variable), 
we include this family in this paper for completeness.    

Following conventional notation, 
a $p$-dimensional random vector $\bY$ 
is said to follow a (multivariate) normal distribution, 
denoted by $\bY \sim N_p(\bmu, \bSigma)$, 
if its density is given by 
\begin{eqnarray}
\phi_p(\by; \bmu, \bSigma) = (2\pi)^{-\frac{p}{2}} |\bSigma|^{-\frac{1}{2}} e^{-\frac{1}{2} (\by-\bmu)^\top \bSigma^{-1} (\by-\bmu)},
\label{eq:N}
\end{eqnarray}
where $\bmu$ is a $p\times 1$ vector of location parameters 
and $\bSigma$ is $p\times p$ positive definite symmetric 
matrix of scale parameters. 
The mean and variance of $\bY$ are $E(\bY)=\bmu$ 
and $\cov(\bY)=\bSigma$, respectively. 
The vector $\bY$ can be expressed as 
a location-scale variant of a standard normal random variable, 
that is,
\begin{equation}
\bY = \bmu + \bSigma^{\frac{1}{2}} \bZ,
\label{eq:Z}
\end{equation}
where $\bZ \sim N_p(\bzero, \bI_p)$, 
$\bzero$ is a vector of zeros, 
and $\bI_p$ is the $p\times p$ identity matrix. 
By `scaling' or  `mixing' $\bY$, 
we mean that $\bmu$ is mixed with $W$ 
and/or $\bSigma$ is weighted by $\sqrt{W}$, 
where $W$ is a positive random variable independent of $\bZ$. 
We  consider each of these cases 
in Sections \ref{VM} to \ref{MVM}. 
By adopting a range of different distributions for $W$, 
a wide variety of non-normal distributions can be constructed.

\section{Variance mixture of normal distributions}
\label{VM}

Variance mixture, or scale mixture, 
of normal (VMN) distributions 
refers to the family of distributions 
generated by scaling the variance matrix $\bSigma$ in (\ref{eq:N}) 
with a (scalar) positive scaling variable $W$. 
More formally, it refers to distributions 
with the following stochastic representation,
\begin{equation}
\bY = \bmu + \sqrt{W} \bSigma^{\frac{1}{2}}\bZ,
\label{eq:VM}
\end{equation}
where $\bZ \sim N(\bzero, \bI_p)$ and $W$ are independent. 
Let the density of $W$ be denoted by $h(w; \btheta)$, 
where $\btheta$ is the vector of parameters associated with $W$. 
It follows that the density is in the form of an integral given by
\begin{equation}
f(\by; \bmu, \bSigma, \btheta) 
= \int_0^\infty \phi_p(\by; \bmu, W\bSigma) \, h(w; \btheta) dw.
\label{eq:VMden}
\end{equation}
A similar expression to (\ref{eq:VMden}) above 
can be given in the case where $W$ has a discrete distribution; 
see, for example, equation (3) of \citet{lm19}. 
As can be observed from (\ref{eq:VM}), 
the family of VMN distributions have constant mean 
but variable scale depending on $W$. 
This allows the VM distributions to have 
lighter or heavier tails than the normal distribution 
and thus are suitable for modelling data 
with tails thickness that deviate from the normal. 
However, this distribution in the unimodal family 
remain symmetric in shape.

\subsection{Properties}
\label{sec:VMprop}

The moments of VMN distributions 
can be readily obtained from (\ref{eq:VM}). 
For example, the first and second moments of $\bY$ 
are given by, respectively, 
$E(\bY)=\bmu$ and $\cov(\bY)= E(W) \bSigma$. 
Further, the moment generating function (mgf) of $\bY$ 
can be expressed as 
\begin{equation}
M_{\bY}(\bt) = e^{\bt^\top\bmu} 
		M_W\left(\frac{1}{2} \bt^\top\bSigma\by\right), 
\label{eq:MVmgf}
\end{equation}
where $M_W(\cdot)$ denotes the mgf of $W$. 

Some nice properties of the normal distribution 
remain valid for VMN distributions, 
including closure under affine transformation, 
marginalization, and conditioning. 
Let $\bY \sim \textit{VMN}_p(\bmu, \bSigma; h(w;\btheta))$ 
denotes $\bY$ having the density (\ref{eq:VM}). 
Let also $\bA$ be a $q\times p$ matrix of full row rank 
and $\ba$ be a $q$-dimensional vector. 
Then the affine transformation $\bA \bY + \ba$ 
still has a VMN distribution with density
\begin{equation}
\bA\bY + \ba \sim 
		\textit{VMN}_q(\bA\bmu + \bb, \bA\bSigma^\top\bA; h(w;\btheta)). 
\label{eq:VMtrans}
\end{equation}
Furthermore, if $\bX \sim VMN_q(\bmu^*, \bSigma^*; h(w;\btheta))$ 
is independent of $\bY$, 
then the linear combination $\bA\bY + \bX$ 
has density given by
\begin{equation}
\bA\bY + \bX \sim 
		\textit{VMN}_q(\bA\bmu + \bmu^*, \bA\bSigma^\top\bA + \bSigma^*; 
				h(w;\btheta)). 
\label{eq:VMcomb}
\end{equation}
Suppose $\bY$ can be partitioned as 
$\bY^\top = (\bY_1^\top, \bY_2^\top)$ 
with respective dimensions $p_1$ and $p_2$ where $p_1+p_2=p$. 
Accordingly, let $\bmu^\top = (\bmu_1^\top, \bmu_2^\top)$ 
and $\bSigma$ be partitioned into four block matrices 
$\bSigma_11$, $\bSigma_{12}$, $\bSigma_{21}$ and $\bSigma_{22}$. 
Then the marginal density of $\bY_1$ 
is $\textit{VMN}_{p_1}(\bmu_1, \bSigma_{11}; h(w;\btheta))$ 
and the conditional density of $\bY_1 | \bY_2=\by_2$ 
is \linebreak $\textit{VMN}_{p_1}(\bmu_{1.2}, \bSigma_{11.2}; h(w;\btheta))$, 
where $\bmu_{1.2}=\bmu_1 + \bSigma_{11}\bSigma_{22}^{-1}(\by_2-\bmu_2)$ 
and $\bSigma_{11.2}=\bSigma_{11}-\bSigma_{12}\bSigma_{22}^{-1}\bSigma_{21}$.

\subsection{Special cases}
\label{sec:VMcases}

The family of VMN distribution encompasses 
many well-known distributions, including the $t$, 
Cauchy, symmetric generalized hyperbolic, 
and logistic distributions. 
The slash, Pearson type VII, contaminated normal, 
and exponential power distributions 
can also be represented as a VMN distribution; 
see also \citet{am74} and \citet{lm19} 
for some other special cases of VMN distributions.

The ($p$-dimensional) $t$-distribution 
can be obtained by 
letting $W\sim IG\left(\frac{\nu}{2}, \frac{\nu}{2}\right)$ 
in (\ref{eq:VM}), where $IG(\cdot)$ denotes 
the inverse gamma distribution and $\nu$ 
is a scalar parameter 
commonly known as the degrees of freedom. 
This tuning parameter regulates 
the thickness of the tails of the $t$-distribution, 
allowing it to model heavier tails than the normal distribution. 
The Cauchy and normal distributions are 
special/limiting cases of the $t$-distribution 
(by letting $\nu=1$ and $\nu\to\infty$, respectively). 

The (symmetric) generalized hyperbolic distribution 
is another important special case of the VMN distribution. 
It arises when $W$ follows 
a generalized inverse Gaussian (GIG) distribution, 
which includes the IG distribution as a special case. 
Thus, the above mentioned $t$-distribution 
and its nested cases are also members of 
the symmetric generalized hyperbolic distribution.

\subsection{Parameter estimation}

From (\ref{eq:VM}), a VMN distribution 
can be expressed in a hierarchical form 
given by $\bY|W=w \sim N_p(\bmu, w\bSigma)$ 
and (with a slight abuse of notation) $W \sim h(w; \btheta)$. 
This facilitates maximum likelihood estimation 
of the model parameters via 
the Expectation-Maximization (EM) algorithm \citep{J034}. 
Technical details can be found in many reports, 
for example, \citet{ls93}.

\section{Mean-mixture of normal distributions}
\label{MM}

Rather than weighting $\bSigma$ with $W$, 
the mean-mixture (or location-mixture) 
of normal (MMN) distribution \citep{njsb19} 
is obtained by mixing $\bmu$ with $W$. 
Note that, in general, $W$ need not be 
a positive random variable in the case MMN distribution. 
More formally, the MMN distribution 
arises from the stochastic expression
\begin{equation}
	\bY = \bmu + W \bdelta + \bSigma^{\frac{1}{2}} \bZ, 
\label{eq:MM}
\end{equation}
where $\bdelta$ is $p\times 1$ vector of shape parameters. 
The density (\ref{eq:MM}) is asymmetric 
if $W$ has an asymmetric distribution. 
In this case, $\bdelta$ may be interpreted as 
a vector of skewness parameters. 
A prominent example is 
the (positively) truncated normal 
or half-normal distribution, that is, 
$W \sim TN(0, 1; \mathbb{R}^+)$. 
This leads to the classical characterization 
of the skew normal (SN) distribution proposed by \citet{J001}.  
It should be noted that while the MMN distribution 
includes the SN distribution as a special case, 
some other commonly used skew-elliptical distributions 
such as the skew $t$-distribution \citep{J006,J019,J012} 
are not MMN distributions.  

From (\ref{eq:MM}), the density of MMN distribution 
can be expressed as 
\begin{equation}
f(\by; \bmu, \bSigma, \bdelta; h(w;\btheta)) 
		= \int_{-\infty}^\infty 
		\phi_p(\by; \bmu + w\delta, \bSigma) \, h(w; \btheta) dw,
\label{eq:MMden}
\end{equation}
where, again, $h(w; \btheta)$ denotes 
the density of $W$ with parameters $\btheta$. 
The notation 
$\bY \sim MMN_p(\bmu, \bSigma, \bdelta; h(w;\btheta))$ 
will be used when $\bY$ has density 
in the form of (\ref{eq:MMden}). 
Similar to the VMN distribution, 
the MMN distribution admits 
a two-level hierarchical representation given by
\begin{equation}
\bY|W=w \sim N_p(\bmu + w\bdelta, \bSigma) 
		\perp W \sim h(w;\btheta). 
\label{eq:MMhier}
\end{equation}

\subsection{Properties}
\label{sec:NNMprop}

It is straightforward to obtain 
the moments for a MMN random variable. 
From the stochastic representation (\ref{eq:MM}), 
it can be seen that first moment of $\bY$ 
is given by $E(\bY)= \bmu+E(W) \bdelta$ if $E(|W|) < \infty$. 
Similarly, the second moment of $\bY$ 
is given by $\cov(\bY)=\bSigma = \var(W) \bdelta\bdelta^\top$, 
provided $E(W^2)$ is finite. 
Further, the mgf of $\bY$ is given by
\begin{equation}
M_{\bY}(\bt) = e^{\bt^\top\bmu+\frac{1}{2} \bt^\top\bSigma\by}
	M_W\left(\bt^\top\bdelta\right).  
\label{eq:MMmgf}
\end{equation} 

The MMN distribution also enjoys nice properties 
such as closure under linear transformation, marginalization, 
and conditioning. 
Let $\bY \sim \textit{MMN}_p(\bmu, \bSigma, \bdelta; h(w;\btheta)$, 
$\bA$ be a $q\times p$ matrix of full row rank, 
and $\ba$ be a $q$-dimensional vector. 
Then the affine transformation $\bA \bY + \ba$ 
still has a NMM distribution with density
\begin{equation}
\bA\bY + \ba \sim 
		\textit{MMN}_q(\bA\bmu + \bb, \bA\bSigma^\top\bA, \bA\bdelta; 
				h(w;\btheta)). 
\label{eq:MMtrans}
\end{equation}
In addition, the linear combination of a MMN 
and a normal random variables is also a MMN random variable. 
If $\bX \sim N_q(\bmu^*, \bSigma^*)$ 
is independent of $\bY$, 
then the linear combination $\bA\bY + \bX$ 
has density given by
\begin{equation}
\bA\bY + \bX \sim 
		\textit{MMN}_q(\bA\bmu + \bmu^*, \bA\bSigma^\top\bA + \bSigma^*, 	 
				\bA\bdelta; h(w;\btheta)). 
\label{eq:MMcomb}
\end{equation}
Concerning the marginal and conditional distributions 
of MMN random variables, let $\bY$, $\bmu$, 
and $\bSigma$ be partitioned 
as in Section \ref{sec:VMprop}. 
Similarly, partition $\blambda$ 
into $\blambda^\top = (\blambda_1^\top, \blambda_2^\top)$. 
Then the marginal density of $\bY_1$ 
is $\textit{MMN}_{p_1}(\bmu_1, \bSigma_{11}, \blambda_1; h(w;\btheta))$ 
and the conditional density of $\bY_1 | \bY_2=\by_2$ 
is $\textit{MMN}_{p_1}(\bmu_{1.2}, \bSigma_{11.2}, \lambda_{1.2}; 
h(w;\btheta))$, 
where $\blambda_{1.2}=\blambda_1 - \bSigma_{12}\bSigma_{22}^{-1}\blambda_2$, 
and $\bmu_{1.2}$ and $\bSigma_{11.2}$ 
are defined in Section \ref{sec:VMprop}.

\subsection{Special cases}

As mentioned previously, 
taking $W \sim TN(0, 1; \mathbb{R}^+)$ 
leads to the classical SN density given by 
\begin{equation}
f(\by; \bmu, \bSigma, \bdelta) = 2 \, \phi_p(\by; \bmu, \bOmega) 
	\Phi_1(\bdelta^\top \bOmega^{-1} (\by-\bmu); 
			0, 1-\bdelta^\top\bOmega^{-1}\bdelta),
\label{eq:SN}
\end{equation} 
where $\bOmega = \bSigma + \bdelta\bdelta^\top$ 
and $\Phi_1(\cdot; \mu, \sigma^2)$ 
denotes the corresponding distribution function of 
$\phi_1(\cdot; \mu, \sigma^2)$. 
When $\bdelta=\bzero$, 
the SN distribution reduces to 
the (multivariate) normal distribution. 

Another special case of the MMN distribution 
were presented in \citet{njsb19}. 
Taking $W$ to have a standard exponential distribution, 
that is $W \sim \exp(1)$, 
leads to the MMN exponential (MMNE) distribution. 
It can be shown that the density is given by
\begin{equation}
f(\by; \bmu, \bSigma, \bdelta) 
	= \frac{\sqrt{2\pi}}{\alpha} e^{\frac{\beta^2}{2}}
	\, \Phi_p(\by; \bmu, \bSigma) \, \Phi_1(\beta),
\label{eq:MMNE}
\end{equation}
where $\alpha^2 = \bdelta^\top \bSigma^{-1}\bdelta$ 
and $\beta=\alpha^{-1}\left[\bdelta^\top\bSigma^{-1}(\by-\bmu)-1\right]$. 
For further details and properties of the MMNE distribution, 
the reader is referred to Section 8.1 in \citet{njsb19}.

\subsection{Parameter estimation}

Utilizing the hierarchical representation (\ref{eq:MMhier}), 
EM algorithm can be implemented 
to provide maximum likelihood estimates of 
the model parameters. 
Although the technical details 
for this EM algorithm are not given 
in \citet{njsb19}, 
it is analogous to the univariate case 
presented in Section 4 of the above reference.

\section{Mean-variance mixture of normal distributions}
\label{MVM}

The mean-variance mixture of normal (MVMN) distribution, 
sometimes called the location-scale mixture of normal 
distribution, is a generalization of 
the VMN distribution described in Section \ref{VM}. 
Compared to (\ref{eq:VM}), the scaling variable $W$ 
is now also mixed with $\bmu$ 
like in the case of the MMN distribution. 
The MVN distribution has 
the following stochastic representation
\begin{equation}
		\bY = \bmu + W \bdelta + \sqrt{W} \bSigma^{\frac{1}{2}} \bZ.
\label{eq:MVMN}
\end{equation}
In this case, both the location and scale of 
the distribution vary with $W$. 
Moreover, $W$ is a positive random variable 
and hence the MVN distribution is asymmetric 
when $\bdelta \neq \bzero$. 
It is important to note that 
while the MVMN distribution reduces to 
the VMN distribution when $\bdelta = \bzero$, 
the MMN distribution described in Section \ref{MM} 
is not a special case of the MVMN distribution.     

Following the definition (\ref{eq:MVMN}), 
the density of $p$-dimensional MVMN distribution 
can be expressed as
\begin{equation}
f(\by; \bmu, \bSigma, \bdelta; h(w;\btheta)) 
		= \int_{0}^\infty 
		\phi_p(\by; \bmu + w\delta, w\bSigma) \, 
		h(w; \btheta) dw. 
\label{eq:MVMden}
\end{equation}
The notation 
$\bY \sim MVMN_p(\bmu, \bSigma, \bdelta; h(w;\btheta))$ 
will be used when $\bY$ has density 
in the form of (\ref{eq:MVMden}). 
Analogous to the VMN and MMN distributions, 
the MMN distribution can be conveniently expressed 
in a hierarchical form given by
\begin{equation}
\bY|W=w \sim N_p(\bmu + w\bdelta, w\bSigma) 
		\perp W \sim h(w;\btheta). 
\label{eq:MVMhier}
\end{equation}

\subsection{Properties}

Some basic properties of the MVMN distribution 
have been studied in \citet{bns82}, among other works. 
The moments of $\bY \sim \textit{MVMN}_p(\bmu, \bSigma, \bdelta; h(w;\btheta))$ 
can be derived directly from (\ref{eq:MVMN}). 
Specifically, the first two moments of $\bY$ 
are given by $E(\bY)-\bmu+E(W)\bdelta$ 
and $\cov(\bY)=\var(W)\bdelta\bdelta^\top+E(W)\bSigma$, 
respectively. Further, the mgf of $\bY$ is given by
\begin{equation}
M_{\bY}(\bt) = e^{\bt^\top\bmu} 
		M_W\left(\bt^\top\bdelta+\frac{1}{2} 
		\bt^\top\bSigma\by\right).  
\label{eq:MVMNmgf}
\end{equation} 

As can be expected, the MVMN distribution 
shares certain nice properties 
with the VMN distribution 
such as closure under linear transformation 
and marginalization. 
Let $\bA$ be a $q\times p$ matrix of full row rank, 
and $\ba$ be a $q$-dimensional vector. 
Then the affine transformation $\bA \bY + \ba$ 
remains a MVMN distribution with density
\begin{equation}
\bA\bY + \ba \sim 
		\textit{MVMN}_q(\bA\bmu + \bb, \bA\bSigma^\top\bA, \bA\bdelta; 
				h(w;\btheta)). 
\label{eq:MVMNtrans}
\end{equation}
Similar to the MMN distribution, 
a linear combination of a MVMN and a normal random variable 
remains a MVMN random variable. 
If $\bX \sim N_q(\bmu^*, \bSigma^*)$ 
is independent of $\bY$, 
then the linear combination $\bA\bY + \bX$ 
has density given by
\begin{equation}
\bA\bY + \bX \sim \textit{MVMN}_q(\bA\bmu + \bmu^*, \bA\bSigma^\top\bA + \bSigma^*, \bA\bdelta; h(w;\btheta)). 
\label{eq:MVMNcomb}
\end{equation}
Marginal distributions and conditional distributions 
of MVMN random variables can also be derived. 
Let $\bY$, $\bmu$, $\bSigma$, and $\bdelta$ 
be partitioned as in Section \ref{sec:NNMprop}. 
Then the marginal density of $\bY_1$ 
is $\textit{MVMN}_{p_1}(\bmu_1, \bSigma_{11}, \blambda_1; h(w;\btheta))$ and the conditional density of 
$\bY_1 | \bY_2=\by_2$ is 
$\textit{MVMN}_{p_1}(\bmu_{1.2}, \bSigma_{11.2}, \lambda_{1.2}; 
h(w;\btheta))$, where $\blambda_{1.2}$, $\bmu_{1.2}$, 
and $\bSigma_{11.2}$ are defined in Section \ref{sec:NNMprop}.

\subsection{Special cases}

Perhaps the most well-known special case of 
the MVMN distribution is 
the generalized hyperbolic (GH) distribution, 
which is widely applied in finance and other fields.
This distribution is obtained by 
letting $W \sim GIG(\psi, \chi, \lambda)$, 
yielding the following density \citep{mfe05},
\begin{equation}
f(\by; \bmu, \bSigma, \bdelta, \psi, \chi, \lambda)
= \frac{\left(\frac{\psi}{\chi}\right)^{\frac{\lambda}{2}} 
		K_{\lambda-\frac{p}{2}}
		\left(\sqrt{(\psi+d_{\bdelta})(\chi+d_{\by})}\right) }
		{(2\pi)^{\frac{p}{2}} |\bSigma|^{\frac{1}{2}} 
		K_\lambda(\chi\psi) 
		e^{\bdelta^\top\bSigma^{-1}(\by-\bmu)}}
		\left(\frac{\chi+d_{\by}}
		{\psi+d_{\bdelta}}\right)^{\frac{\lambda}{2}-\frac{p}{4}},
\label{eq:GH}
\end{equation}
where $d_{\bdelta} = \bdelta^\top\bSigma^{-1}\bdelta$, 
$d_{\by}=(\by-\bmu)^\top \bSigma^{-1} (\by-\bmu)$, 
and $K_\lambda(\cdot)$ denotes 
the modified Bessel function of the third kind 
with index $\lambda$. 
The GH distribution, as the name suggests, 
contains the symmetric GH distribution 
mentioned in Section \ref{sec:VMcases} 
and an asymmetric version of some of its members. 
However, it cannot obtain the SN distribution 
as a special/limiting case. 
Other noteworthy special cases of the GH distribution 
include the normal inverse Gaussian, 
variance gamma, and asymmetric Laplace distributions. 
The GH distribution and its properties 
have been well studied in the literature; 
see, for example, \citet{T004} and \citet{dy18}.  

Two other less well-known MVMN distributions 
were recently considered by \citet{pjr15} and \citet{naj18}. 
The former presented 
a MVN of Birnbaum-Saunders (MVNBS) distribution, 
where $W$ has a Birnbaum-Saunders distribution with shape parameter $\alpha$ and scale parameter $1$. 
In the second reference, the authors assumed 
$W$ follows a Lindley distribution, 
which is a mixture of $\exp(\alpha)$ 
and $\mbox{gamma}(2, \alpha)$ distributions. 
This leads to the so-called MVN Lindley (MVNL) distribution.

\subsection{Parameter estimation}

The EM algorithm can be employed 
to estimate the parameters of the MVMN distribution. 
For special cases of MVMN distribution 
such as the GH, MVNBS, and MVNL distributions, 
explicit expressions for the implementation 
of the EM algorithm can be found in 
\cite{J407,pjr15}, and \citet{naj18}, respectively.

\section{Conclusions}
\label{concl}

A concise description of three generalizations 
of the (multivariate) normal distribution 
has been presented. 
These families of flexible distributions 
arise by mixing the mean and/or 
weighting the variance matrix 
of a normal random variable. 
Two of these families,
namely the variance mixture (VMN) 
and mean-variance mixture of normal (MVMN) distributions 
have a relatively long history in the literature, 
whereas the third family 
(mean-mixture of normal (MMN) distribution) 
were introduced more recently. 
Each of these families has their own merits and limits. 
We have presented their basic properties, 
some important special/limiting cases, 
and references for parameter estimation procedures. 
Some further versions and/or generalizations 
of MVMN would be of interest for future investigation; 
for example, 
a scale mixture of MMN distributions 
(as suggested by \cite{njsb19}) 
and a MVMN distribution 
where different mixing variables can be used 
for the mean and variance of the normal random variable.



\begin{thebibliography}{21}
\providecommand{\natexlab}[1]{#1}
\providecommand{\url}[1]{\texttt{#1}}
\providecommand{\urlprefix}{URL }

\bibitem[{Adcock and Azzalini(2020)}]{aa20}
Adcock, C. and Azzalini, A. (2020).
\newblock A selective overview of skew-elliptical and related distributions and
  of their applications.
\newblock \emph{symmetry} \textbf{12}, 118.

\bibitem[{Andrews and Mallows(1974)}]{am74}
Andrews, D.F. and Mallows, C.L. (1974).
\newblock Scale mixtures of normal distributions.
\newblock \emph{Journal of the Royal Statistical Society, Series B}
  \textbf{36}, 99--102.

\bibitem[{Arellano-Valle and Azzalini(2006)}]{J017}
Arellano-Valle, R.B. and Azzalini, A. (2006).
\newblock On the unification of families of skew-normal distributions.
\newblock \emph{Scandinavian Journal of Statistics} \textbf{33}, 561--574.

\bibitem[{Azzalini(2005)}]{J021}
Azzalini, A. (2005).
\newblock The skew-normal distribution and related multivariate families.
\newblock \emph{Scandinavian Journal of Statistics} \textbf{32}, 159--188.

\bibitem[{Azzalini and Capitanio(2003)}]{J006}
Azzalini, A. and Capitanio, A. (2003).
\newblock Distributions generated by perturbation of symmetry with emphasis on
  a multivariate skew $t$-distribution.
\newblock \emph{Journal of the Royal Statistical Society B} \textbf{65},
  367--389.

\bibitem[{Azzalini and Capitanio(2014)}]{B025}
Azzalini, A. and Capitanio, A. (2014).
\newblock \emph{The Skew-Normal and Related Families}.
\newblock Cambridge: Cambridge University Press.

\bibitem[{Azzalini and {\uppercase{D}alla Valle}(1996)}]{J001}
Azzalini, A. and {\uppercase{D}alla Valle}, A. (1996).
\newblock The multivariate skew-normal distribution.
\newblock \emph{Biometrika} \textbf{83}, 715--726.

\bibitem[{Barndorff-Nielsen et~al.(1982)Barndorff-Nielsen, Kent, and
  S{\o}rensen}]{bns82}
Barndorff-Nielsen, O., Kent, J., and S{\o}rensen, M. (1982).
\newblock Normal variance-mean mixtures and z distributions.
\newblock \emph{International Statistical Review} \textbf{50}, 145--159.

\bibitem[{Branco and Dey(2001)}]{J012}
Branco, M.D. and Dey, D.K. (2001).
\newblock A general class of multivariate skew-elliptical distributions.
\newblock \emph{Journal of Multivariate Analysis} \textbf{79}, 99--113.

\bibitem[{Browne and McNicholas(2015)}]{J407}
Browne, R.P. and McNicholas, P.D. (2015).
\newblock A mixture of generalized hyperbolic distributions.
\newblock \emph{The Canadian Journal of Statistics} \textbf{43}, 176--198.

\bibitem[{Dempster et~al.(1977)Dempster, Laird, and Rubin}]{J034}
Dempster, A.P., Laird, N.M., and Rubin, D.B. (1977).
\newblock Maximum likelihood from incomplete data via the {EM} algorithm.
\newblock \emph{Journal of Royal Statistical Society B} \textbf{39}, 1--38.

\bibitem[{Deng and Yao(2018)}]{dy18}
Deng, X. and Yao, J. (2018).
\newblock On the property of multivariate generalized hyperbolic distribution
  and the stein-type inequality.
\newblock \emph{Communications in Statistics - Theory and Methods} \textbf{47},
  5346--5356.

\bibitem[{Genton(2004)}]{B001}
Genton, M.G. (Ed.). (2004).
\newblock \emph{Skew-Elliptical Distributions and Their Applications: A Journey
  Beyond Normality}.
\newblock Boca Raton, Florida: Chapman \& Hall, CRC.

\bibitem[{Gupta(2003)}]{J019}
Gupta, A.K. (2003).
\newblock Multivariate skew-$t$ distribution.
\newblock \emph{Statistics} \textbf{37}, 359--363.

\bibitem[{Iversen(1999)}]{T004}
Iversen, D. (1999).
\newblock \emph{The Generalized Hyperbolic Model:Estimation, Financial
  Derivatives, and Risk Measures}.
\newblock Master's thesis, Albert-Ludwigs-Universit{\"a}t Freiburg.

\bibitem[{Lange and Sinsheimer(1993)}]{ls93}
Lange, K. and Sinsheimer, J.S. (1993).
\newblock Normal/independent distributions and their applications in robust
  regression.
\newblock \emph{Journal of Computational and Graphical Statistics} \textbf{2},
  175--198.

\bibitem[{Lee and McLachlan(2019)}]{lm19}
Lee, S. and McLachlan, G. (2019).
\newblock Scale mixture distribution.
\newblock \emph{Wiley Stats Ref: Statistics Reference Online (WSR)} 08201.

\bibitem[{McNeil et~al.(2005)McNeil, Frey, and Embrechts}]{mfe05}
McNeil, A.J., Frey, R., and Embrechts, P. (2005).
\newblock \emph{Quantitative Risk Management: Concepts, Techniques and Tools}.
\newblock New Jersey, US: Princeton University Press.

\bibitem[{Naderi et~al.(2018)Naderi, Arabpour, and Jamalizadeh}]{naj18}
Naderi, M., Arabpour, A., and Jamalizadeh, A. (2018).
\newblock Multivariate normal mean-variance mixture distribution based on
  {L}indley distribution.
\newblock \emph{Communications in Statistics-Simulation and Computation}
  \textbf{47}, 1179--1192.

\bibitem[{Negarestani et~al.(2019)Negarestani, Jamalizadeh, Shafiei, and
  Balakrishnan}]{njsb19}
Negarestani, H., Jamalizadeh, A., Shafiei, S., and Balakrishnan, N. (2019).
\newblock Mean mixtures of normal distributions: properties, inference and
  application.
\newblock \emph{Metrika} \textbf{82}, 501--528.

\bibitem[{Pourmousa et~al.(2015)Pourmousa, Jamalizadeh, and Rezapour}]{pjr15}
Pourmousa, R., Jamalizadeh, A., and Rezapour, M. (2015).
\newblock Multivariate normal mean variance mixture distribution based on
  {B}irnbaum {S}aunders distribution.
\newblock \emph{Journal of Statistical Computation and Simulation} \textbf{85},
  2736--2749.

\end{thebibliography}
\end{document}